\documentclass[10pt, a4paper, twoside]{amsart}

\usepackage{latexsym}
\usepackage[all]{xy}
\usepackage{color}
\newtheorem{teorema}{Theorem}[section]
\newtheorem{definicion}[teorema]{Definition}
\include{amslatex}
\newtheorem{proposicion}[teorema]{Proposition}

\newtheorem{corolario}[teorema]{Corollary}

\newtheorem{comentario}[teorema]{Remark}

\newtheorem{ejemplo}[teorema]{Example}

\newtheorem*{theorem-A}{Theorem A}
\newtheorem*{theorem-B}{Theorem B}
\newtheorem*{theorem-C}{Theorem [Akdar-Zadeh]}

\numberwithin{equation}{section}

\title[Average structures]{Average structures associated to a Finsler space}

\author[Gallego Torrom\'e]{Ricardo Gallego Torrom\'e}
\email{rigato39@gmail.com}
\address{Departamento de Matem\'atica\hfill\break\indent Universidade Federal de S\~ao Carlos\hfill\break\indent Brazil}
\date{March 2sd, 2017}
\begin{document}

\maketitle

\begin{abstract}
Given a Finsler space $(M,F)$ on a manifold $M$, the averaging method associates to  {\it Finslerian geometric objects} {\it affine geometric objects} living on $M$. In particular, a Riemannian metric is associated to the fundamental tensor $g$ and an affine, torsion free connection is associated to the Chern-Rund connection. As an illustration of the technique, a generalization of the Gauss-Bonnet theorem to Berwald surfaces using the average metric is presented. The parallel transport and curvature endomorphisms of the average connection are obtained. The holonomy group for a Berwald space is discussed.  New affine, local isometric invariants of the original Finsler metric. The heredity of the property of symmetric space from the Finsler space to the average Riemannian metric is proved.
\end{abstract}

\begin{section}{Introduction}

It is  notable that many fundamental results can be generalized from the Riemannian to the Finslerian category by properly adapting the proofs, in some cases in a rather straightforward way \cite{BCS}. The present paper is motivated by the idea of interpreting this general phenomenon from the perspective of equivalence classes in the Finsler category such that each class contains an affine or Riemannian representative. The reason for such expectation is that, if a given proposition can be casted in terms of notions defined in the coset space of equivalent classes, then it can be investigated using affine or Riemannian methods.

The method that we use to define the equivalence relations is by  {\it averaging Finslerian objects}. These operations transform geometric objects living in the tangent space to geometric objects living in the base manifold $M$. The average of the fundamental tensor appeared first in \cite{Csava2005}, in connection with emergent quantum mechanics \cite{Ricardo 2006} and in a series of pre-prints of the author (RGT), of which the present work is the latest version. Preliminary versions of the general theory of averaging were also discussed by the author, for instance in \cite{Master Ricardo}. A revision of the theory is presented in this paper. We will show that we can consider the average of the Finslerian metric (in the form of the fundamental tensor), connection, curvatures or the average differential operators. Two geometric operators/structures will be  equivalent if they have the same average. It turns out that results stated  in terms of the defining properties of the coset space and that can be proved using Riemannian or affine methods automatically upgrade to the Finsler category. In particular, the method is very well suited for two specific types of problems. The fist is in applications to Berwald spaces. The second in the application to isometry properties of the Finsler metric and related notions.

In this paper we discuss the general method of average and its application to several geometric operators and structures. These applications are illustrated by several examples.
This paper is organized as follows. In {\it section 2}, several standard notions of Finsler geometry are introduced.
In {\it section 3} we explain the measure used for the averaging operation. The averaging can be casted in terms of
integration  along vertical fibers of  pull-back bundles $\hat{\pi}^* { T^{(p,q)}M}$ over the  manifold $TM\setminus
\{ 0\} $. The averaging associates a  Riemannian structure $(M, h)$ to the initial Finsler space $(M,F)$ by performing an average of the components of the fundamental tensor $g$ on the indicatrix ${I}_x $ at each point $ x\in M$.
In {\it section 4} we discuss the average connection $\langle \nabla \rangle$ associated to  a linear connection  $\nabla$ on $\hat{\pi}^*TM$. We show that for Berwald spaces, the Riemann curvature tensor of the average metric $h$ can be obtained by averaging the corresponding tensor of the Chern-Rund connection of $g$.  The average metric is used to prove several generalizations of Riemannian results to Berwald spaces. In particular, we have considered a  generalization of the Gauss-Bonnet theorem for Berwald surfaces as an example of application of the method.
In {\it section 5} the parallel transport operators and curvatures of the average connection of $\langle \nabla \rangle$ are discussed in terms of the corresponding parallel transport and curvatures of the linear connection  $\nabla$ on $\hat{\pi}^*TM$. The metrizability of the holonomy of Berwald spaces is discussed. In {\it section 6}, isometries are discussed and the proof of the  Mayer-Steenrod theorem sketched from the point of view of the averaging method. Symmetric spaces are discussed too from this new framework. The average of the curvature tensor of the connection $\nabla$ provides new affine isometric invariants of the Finsler metric.
\end{section}

\begin{section}{Basic notions of Finsler geometry}
\subsection*{Notation}
 $(U,{\bf x})$ will denote a local coordinate chart of the $n$-dimensional manifold $M$, where a point $x\in U$ has local coordinates $(x^{1},...,x^{n})$ and $U\subset M$ is an open set.
$TM$ is the tangent bundle of $M$. The slit
tangent bundle $\hat{\pi}:N\to M$ is the bundle over $M$  with
$N=TM\setminus \{0\}$.
 Fixed a local chart $(U,{\bf x})$ on the manifold $M$, a point $x\in\, U$ will have coordinates $(x^{1},...,x^{n})$ and a tangent vector $y=y^{i}\frac{\partial}{\partial x^i}\in T_x 
M$ at $x\in M$ is determined by its components $y=(y^1 ,...,y^n )$ respect to the basis $\{\frac{\partial}{\partial x^i}\}^n_{i=1}$of $T_xM$. Note that we will use Einstein's convention 
for up and down equal indices, if anything else is not
directly stated. The set of sections
of a bundle $\mathcal{E}$ is denoted by $\Gamma \,\mathcal{E}$, except for differential forms that we follow the usual notation. Each
local chart $(U,{\bf x})$ on  $M$
induces a local chart on $TM$ denoted by
$(TU,{\bf x},{\bf y})$ such that a point $u\,\in TU$ with $\hat{\pi}(u)=x$ and corresponding to the tangent vector $y=y^{i}\frac{{\partial}}{{\partial}
x^{i}}\in T_xM$ has local {\it natural coordinates}
$(x^1,...,x^n ,y^1,...,y^n)$.
\begin{definicion}
\label{Finslerstructute}
A Finsler space on the manifold $M$  is a pair $(M,F)$ where $F$ is a
non-negative, real function  $F:M\to [0,\infty [$
such that
\begin{itemize}
\item It is smooth in the slit tangent bundle $N$,

 \item Positive homogeneity holds: $F(x,{\lambda}\,y)=\,\lambda \,F(x,y)$ for every $\lambda \in\,[0,+\infty [$,

 \item Strong convexity holds:
the Hessian matrix
 \begin{align}
g_{ij}(x,y): =\,\frac{1}{2}\,\frac{{\partial}^2 F^2
(x,y)}{{\partial}y^i {\partial}y^j },\quad i,j=\,1,...,n
\label{fundamentaltensorcoefficients}
\end{align}
is positive definite on $N$.
\end{itemize}
\label{definicionfinslerstructure}
\end{definicion}
 The minimal regularity requirement for the Finsler function $F$ is to be a  $\mathcal{C}^4$-smooth function on $N$.
However, if the second Bianchi identities are required, then it is necessary for  $F$ to be at least $\mathcal{C}^5$-smooth on $N$. The matrix $(g)_{ij}:=g_{ij}(x,y)$ is the matrix components of the fundamental tensor $g$ at the point $u=(x,y)\in\,N$.
\begin{definicion}
Let $(M, F)$ be a Finsler space and $(TU,{\bf x},{\bf y}) $ a
local chart induced on $N$ from the coordinate system $(U,{\bf x})$ of 
$M$. The components of the Cartan tensor are defined by the
collection of functions
\begin{align}
{ A}_{ijk}=\frac{F}{2} \frac{\partial g_{ij}}{\partial
y^{k}},\quad i,j,k=1,...,n.
\label{cartancoefficients}
\end{align}
\end{definicion}
The components  $\{A_{ijk},\,i,j,k=1,...,n\}$ are homogeneous functions of degree zero in the coordinates $(y^1,...,y^n )$ and totally symmetric under permutation of the indices $i,j,k$. The condition
\begin{align*}
A_{ijk}(x,y)=0,\quad\forall\, y\in T_xM,\,x\in\,M
\end{align*}
  characterizes Riemannian geometry among the general class of Finsler geometries.
\begin{definicion}
Let $(M,F)$ be a Finsler space. The indicatrix over the point $x\in M$ is the submanifold $I_x\hookrightarrow T_xM$
\begin{align*}
I_x :=\{ y\in T_x M\mid F(x,y)=1\}
\end{align*}
\label{definitionofindicatrix}
\end{definicion}
The indicatrix $I_x$ is a compact, strictly convex submanifold of $T_x M$ \cite{BCS, Rund}.
 We denote by $\mathcal{I}$ the fibered manifold $\pi_{\mathcal{I}}:\mathcal{I}\to M$ with $\pi^{-1}_{\mathcal{I}}(x)=\,I_x$ and the base manifold is $M$.

\subsection{Definition of a  non-linear connection of $N$}
Let us consider the slit bundle $\hat{\pi}:N\to M$.
\begin{definicion}
A non-linear connection of $N$ is a distribution $\mathcal{H}\subset \,TN$  supplementary to the canonical vertical distribution $\mathcal{V}=\,\ker \, d\hat{\pi}$.
\label{nonlinearconnection}
\end{definicion}
Given a Finsler space $(M,F)$, there is defined a non-linear connection in the manifold $N$.
In a local natural coordinate chart $(TU,{\bf x},{\bf y}) $ of $N$ the collection of local sections
\begin{align*}
\left\{ \frac{\partial}{\partial 
y^{1}}|_u ,..., \frac{\partial}{\partial y^{n}}|_u ,\,u\in\hat{\pi}^{-1}(x),\, x\in U \right\}
\end{align*}
determines a local frame for the vertical distribution $\mathcal{V}$. To obtain a supplementary distribution $\mathcal{H}$ we use a standard local construction \cite{BCS}. First, let us introduce the non-linear
connection coefficients ${N^{i}_{j}}$  by the expression
\begin{align*}
\frac{N^{i}_{j}}{F}={\gamma}^{i}_{jk}\frac{y^{k}}{F}-A^{i}_{jk}
{\gamma}^{k}_{rs}\frac{y^{r}}{F}\frac{y^{s}}{F},\quad
i,j,k,r,s=1,...,n,
\end{align*}
where the {\it formal second kind Christoffel's symbols}
${\gamma}^{i}_{jk}(x,y)$ are defined by
\begin{align*}
 {\gamma}^{i}_{jk}=\frac{1}{2}g^{is}(\frac{\partial g_{sj}}{\partial
x^{k}}-\frac{\partial g_{jk}}{\partial x^{s}}+\frac{\partial
g_{sk}}{\partial x^{j}}),\quad i,j,k=1,...,n
\end{align*}
and also $A^i _{jk}:=g^{il}A_{ljk}$ and $g^{il}g_{lj}=\delta ^i _j .$  A tangent basis for $T_u N$ is determined
by the  vectors
\begin{align}
\left\{ \frac{{\delta}}{{\delta} x^{1}}|_u ,...,\frac{{\delta}}{{\delta} x^{n}} |_u, 
F\frac{\partial}{\partial y^{1}} |_u,...,F\frac{\partial}{\partial
y^{n}} |_u\right\},\,
 \frac{{\delta}}{{\delta} x^{j}}|_u =\frac{\partial}{\partial
x^{j}}|_u -N^{i}_{j}\frac{\partial}{\partial y^{i}}|_u ,\,
i,j=1,...,n.
\label{localframeTN}
\end{align}
The collection of local sections
\begin{align*}
\left\{ \frac{{\delta}}{{\delta} x^{1}}|_u ,...,\frac{{\delta}}{{\delta} x^{n}}|_u,\, u\in \hat\pi^{-1}(x),\, x\in U \right\}
\end{align*}
  determines a local frame for the {\it horizontal distribution} $\mathcal{H}$ \cite{Bao, BCS}. Given $\tilde{X}\in\,\Gamma TN$, the horizontal component is denoted by $H(\tilde{X})$ and the vertical component by $V(\tilde{X})$.
The horizontal lift of tangent vectors is defined by the homomorphism
\begin{align}
\iota_u:\,T_xM \to  \, T_u N, \quad
  X=X^i\frac{{\partial}}{{\partial}x^i}|_x \mapsto
\iota_u(X)=X^i\frac{{\delta}}{{\delta}x^i}|_{u},
\label{horizontalliftofvectorfields}
\end{align}
for  fixed $\, u\in \,\pi^{-1}(x).$
For these local horizontal sections (and therefore for any local horizontal section) the relation
\begin{align}
\frac{{\delta}}{{\delta} x^{1}}|_u\,\cdot F=0
\label{horizontalderivativeofF}
\end{align}
holds good.

The dual basis associated to the local tangent basis  \eqref{localframeTN} of $T_uN$ is a local basis of the dual vector space ${ T^{*}_u N}$,
\begin{align}
\left\{ dx^{1}|_u ,...,dx^{n}|_u, \frac{{\delta}y^{1}}{F}|_u,...,
\frac{{\delta}y^{n}}{F}|_u \right\},\quad
\frac{{\delta}y^{i}}{F}|_u =\frac{1}{F}(dy^{i}+N^{i}_{j}dx^{j})|_u, \quad i,j=1,...,n.
\end{align}

\begin{definicion}
Let $(M,F)$ be a Finsler space. The fundamental and the Cartan tensors are defined 
in the natural local coordinate system $(TU,{\bf x},{\bf y})$ by the
equations
\begin{itemize}
\item The fundamental tensor is
\begin{align}
g(x,y):=\frac{1}{2}\frac{{\partial}^2 F^2(x,y) }{{\partial}y^i
{\partial}y^j }\, dx^i \otimes dx^j ,
\label{fundamentaltensor}
\end{align}
\item The Cartan tensor is
\begin{align}
A(x,y):=\frac{F}{2} \frac{\partial g_{ij}}{\partial y^{k}}\,
\frac{{\delta}y^i}{F} \otimes dx^j \otimes dx^k=A_{ijk}\,
\frac{{\delta}y^i}{F} \otimes dx^j \otimes dx^k .
\label{cartantensor}
\end{align}
\end{itemize}
\end{definicion}
Note that our definition of the Cartan's tensor \eqref{cartantensor} differs from the usual definition \cite{BCS}, which is a symmetric tensor $A=\,A_{ijk}\,dx^i \otimes dx^j \otimes dx^k $. Such difference is only formal and does not spoil the properties of the Cartan's tensor in any practical situation.

\subsection{The Chern-Rund connection}
Let us consider the  product ${ N}\times { TM} $ and the canonical projections
\begin{align*}
\hat\pi _1:{ N}\times { TM}\to { N},\quad
(u,\xi)\mapsto u,\quad\quad
\hat\pi _2 :{ N}\times { TM}\to {TM},\quad
(u,\xi)\mapsto \xi.
\end{align*}
The pull-back $\hat{\pi}_1:\hat{\pi}^* { TM}\to { N}$ of the tangent bundle
$\pi:{TM}\to {M}$ by the projection $\hat{\pi}:{N}\to { M}$
 is the maximal sub-manifold of the product
${N}\times {TM} $ such that  $\pi \circ\hat\pi
_2(u,\xi)=\hat{\pi}\circ \hat\pi_1(u,\xi)$ holds. It follows that the diagram
\begin{displaymath}
\xymatrix{\hat{\pi}^*{ TM} \ar[d]_{\pi_1} \ar[r]^{\pi_2} &
{TM} \ar[d]^{\pi}\\
{ N} \ar[r]^{\hat{\pi}} & { M}}
\end{displaymath}
 commutes.
Here $\pi_1:\hat{\pi}^* { TM}\to { N}$  and $\pi_2:\hat{\pi}^* { TM}\to { TM}$ are the restrictions of the natural projections $\hat{\pi}_1:{ N}\times { TM}\to { N}$ and $\hat{\pi}_2 :{ N}\times { TM}\to {TM}$ to $\hat{\pi}^* { TM}$.
$\pi_1:\hat{\pi}^* { TM}\to { N}$ is a real vector bundle,
with fiber over $u=(x,Z)\in { N}$ isomorphic to  ${ T}_x{ M}$.
The fiber dimension of $\hat{\pi}^*{TM}$ is equal to $n=dim({ M})$, while the dimension of the base manifold $N$ is $2n$.

Each tangent field $Z\in { T}_x{ M}$ can be pulled-back $\hat{\pi}^*Z$ uniquely by the conditions
 \begin{itemize}
 \item $\pi_1(\hat{\pi}^*Z)=(x,Z)\in {N}_x$,

 \item $\pi_2(\hat{\pi}^* Z)=Z$.
 \end{itemize}
 These conditions extend pointwise to vector fields.
 The pull-back of a smooth function $f\in \mathcal{F}({M})$ is the smooth function $\hat{\pi}^*f\in \mathcal{F}(N)$ such that $\hat{\pi}^*f(u)=f(\hat{\pi}(u))$ for every $u\in { N}$.
 Analogously, a pull-back tensor bundle $\hat{\pi}^*T^{(p,q)}M$ can be obtained from each tensor bundle $T^{(p,q)}M$ over ${ M}$.

Let us consider a linear connection on $\hat{\pi}^*TM$.
The associated covariant derivative is the operator
\begin{align*}
\nabla: \,\Gamma\,\hat{\pi}^*TM\times \,\Gamma \,{T}N\to\, \Gamma\,\hat{\pi}^*TM
\end{align*}
such that
\begin{itemize}
\item For every $\tilde{X}\in \Gamma TN$, $S_1, S_2\in \Gamma\,\hat{\pi}^*TM$ and $\varphi\in\mathcal{F}(M)$ it holds that
\begin{align}
\nabla_{\tilde{X}} (\hat{\pi}^*\varphi S_1+S_2)=\,(\tilde{X}\cdot \hat{\pi}^* \phi)S_1+\,\varphi\nabla_{\tilde{X}} S_1+ \,\nabla_{\tilde{X}} S_2.
\label{leibnitzrule}
\end{align}

\item For every $\tilde{X}_1, \tilde{X}_2\in \Gamma\, TN$, $S\in \Gamma\,\hat{\pi}^*TM$ and $\tilde{\varphi}\in\mathcal{F}(N)$ it holds that
\begin{align}
\nabla_{\tilde{\varphi} \tilde{X}_1+\tilde{X}_2}S=\,\tilde{\varphi}\nabla_{\tilde{X}_1}S+\,\nabla_{\tilde{X}_2}S.
\label{Xlinearity}
\end{align}
\end{itemize}
The Chern-Rund connection $^{ch}\nabla$ is a linear connection on  $\hat{\pi}^*TM$ characterized by the following \cite{BCS, Rund}
\begin{teorema}
Let $(M,F)$ be a Finsler space. The vector bundle $\hat{\pi}^{*}TM$ admits a 
unique linear connection characterized by the collection of connection 1-forms $\{ {\,^{ch}\omega}^i _j,\,\, 
i,j=1,...,n \} $ such that the following structure equations hold:
\begin{itemize}
\item ``Torsion free" condition,
\begin{align}
 d(dx^{i})-dx^{j}\wedge
\,^{ch}w^{i}_{j}=0,\quad i,j=1,...,n.
\label{torsionfreecondition}
\end{align}

\item Almost g-compatibility condition,
\begin{align}
dg_{ij}-g_{kj}\,^{ch}w^{k}_{i}-g_{ik}\,^{ch}w^{k}_{j}=2A_{ijk}\frac{{\delta}y^{k}}{F},\quad
i,j,k=1,...,n.
\label{almostgcompatibility}
\end{align}
\end{itemize}
\label{theoremfromchern}
\end{teorema}
The torsion free condition is equivalent to the absence of terms
containing $\delta y^{i}$ in the connection 1-forms ${\,^{ch}\omega}^i _j $ and
also implies the symmetry of the connection coefficients $\Gamma^i_{jk}(x,y)$,
\begin{align}
^{ch}\omega^{i}_{j}(x,y)={\Gamma}^{i}_{jk}(x,y)\, dx^{k},\quad
{\Gamma}^{i}_{jk}(x,y)={\Gamma}^{i}_{kj}(x,y),\quad i,j,k=1,...,n.
\label{torsionpropertiesofchernconnection}
\end{align}

 We can characterize the Chern-Rund connection by means of the associated covariant derivative operator $^{ch}\nabla_{\tilde{X}}$, with $\tilde{X}\in\,T_u N$.  The {\it generalized torsion tensor} is given by the expression
\begin{align}
T_{\,^{ch}\nabla}:  &\Gamma\, T M \times \Gamma\, TM\to\Gamma\, TM\nonumber \\
 & (X,Y)\mapsto \pi_2\big({^{ch}\nabla}_{\tilde{X}} \hat{\pi}^* Y\big)-\pi_2\big({^{ch}\nabla}_{\hat{Y}}\hat{\pi}^*X\big)-[X,Y],
\label{generalizedtorsiontensor0}
\end{align}
where $\tilde{X}$ and $\tilde{Y}$ are the horizontal lifts of the restrictions $X(x),Y(x)$ at a given point $u\in\,\hat{\pi}^{-1}(x)$.
Then the following corollaries are direct consequences of {\it Theorem} \ref{theoremfromchern}.
\begin{corolario}
Let $(M, F)$ be a Finsler space. The torsion free condition \eqref{torsionfreecondition} is equivalent to the following conditions:
\begin{enumerate}
\item For any $\tilde{X} \in { TN}$ and $Y\in M$, the following relation holds,
\begin{align}
{^{ch}\nabla}_{{V(\tilde{X})}} \hat{\pi}^* Y=0.
\end{align}
\item  Let us consider $X,Y\in \Gamma TM$. Then the following relation holds,
\begin{align}
T_{\,^{ch}\nabla}(X,Y)=0.
\label{generalizedtorsiontensor}
\end{align}
\end{enumerate}
\end{corolario}
\begin{proof}Let $\{\frac{\partial }{\partial x^i}\}^n_{i=1}$ be a local frame on $U$ and $X=e_i=\,\frac{\partial}{\partial x^i}$, $Y=e_j=\frac{\partial}{\partial x^j}.$ Then we have
\begin{align*}
{^{ch}\nabla}_{{V(\tilde{X})}} \hat{\pi}^* Y  :=\hat{\pi}^* e_k\,w^k_j(V(e_i))=\hat{\pi}^*
e_k \Gamma^k_{aj}dx^a(V(e_i)) =\hat{\pi}^* e_k
\Gamma^k_{aj}dx^a(\frac{\partial}{\partial y^i})=0.
\end{align*}
This condition is extended by  linearity and by the Leibnitz rule to arbitrary sections $X,\, Y\in \Gamma\,TM$.

To prove the relation \eqref{generalizedtorsiontensor}, let us consider the torsion condition in local
coordinates where the local frame $\{e_j\}$ on $\Gamma TM$ commutes, $[e_i,e_j]=0$. Then from the symmetry in the connection
coefficients \eqref{torsionpropertiesofchernconnection} one obtains
\begin{align*}
\pi_2({^{ch}\nabla}_{\tilde{e_i}} \hat{\pi}^* e_j)-\,\pi_2({^{ch}\nabla}_{\tilde{e_j}}\hat{\pi}^*
e_i)-[e_i,e_j]& =\,\pi_2({^{ch}\nabla}_{\tilde{e}_i} \hat{\pi}^*
e_j)-\,\pi_2({^{ch}\nabla}_{\tilde{e}_j}\hat{\pi}^*e_i)\\
& =\,(\Gamma^k_{ij}-\Gamma^k_{ji})\,\hat{\pi}^*e_k=0.
\end{align*}
This relation is extended by linearity to arbitrary vectors $X,Y\in \,\Gamma TM$.
\end{proof}
\begin{corolario}
Let $(M, F)$ be a Finsler space and $\tilde{X}\in \,\Gamma\,{T N}$. The almost $g$-compatibility condition \eqref{almostgcompatibility} is equivalent to the conditions
\begin{enumerate}
\item ${^{ch}\nabla}$ is metric compatible in the horizontal directions,
\begin{align}
 {^{ch}\nabla}_{H(\tilde{X})} g=0.
 \label{verticalmetriccondition}
\end{align}
\item ${^{ch}\nabla}$ is almost-metric compatible in the vertical directions in the sense that
\begin{align}
{^{ch}\nabla}_{V(\tilde{X})}g=\,2\,A(\tilde{X},\cdot,\cdot)
\label{horizontalmetricondition}
\end{align}
holds good.
\end{enumerate}
\end{corolario}
\begin{proof} Using local natural coordinates and reading from the expression \eqref{almostgcompatibility}, it follows that
\begin{align*}
^{ch}\nabla g=(dg_{ij}-g_{kj}\,^{ch}w^{k}_{i}-g_{ik}\,^{ch}w^{k}_{j})\,\hat{\pi}^*e^i\otimes
\hat{\pi}^*e^j\,=\, 2\,A_{ijk}\frac{{\delta}y^{k}}{F}\otimes \hat{\pi}^*
e^i\otimes\hat{\pi}^* e^j,
\end{align*}
By the definition of covariant derivative along an horizontal direction and
since $\, 2A_{ijk}\frac{{\delta}y^{k}}{F}$ is vertical, the relation \eqref{verticalmetriccondition} follows. For the covariant derivative along the vertical component $V(\tilde{X})$,
\begin{align*}
^{ch}\nabla_{V(\tilde{X})} \,g:=\, 2\,A_{ijk}\frac{{\delta}y^{k}}{F}\cdot V(\tilde{X})\,\left(\hat{\pi}^*
e^i\otimes \hat{\pi}^* e^j\right),\,\,\forall \tilde{X} \in \Gamma { TN},
\end{align*}
from where follows the condition \eqref{horizontalmetricondition}.
\end{proof}
The curvature $2$-forms associated with a linear connection $\nabla$ on $\hat{\pi}^*TM$  are
\begin{equation}
{\Omega}^{i}_{j}:=dw^{i}_{j}-w^{k}_{j}\wedge w^{i}_{k},\quad
i,j,k=1,...,n.
\end{equation}
In local coordinates, the curvature endomorphisms are decomposed as
\begin{equation}
{\Omega}^{i}_{j}=\frac{1}{2} R^{i}_{jkl}\,dx^{k}\wedge dx^{l}+
P^{i}_{jkl}\,dx^{k}\wedge
\frac{{\delta}y^{l}}{F}+\frac{1}{2}Q^{i}_{jkl}\,\frac{{\delta}y^{k}}{F}\wedge
\frac{{\delta}y^{l}}{F}.
\end{equation}
The quantities $R^i _{jkl}$, $P^i _{jkl}$ and $Q^i _{jkl}$ are the hh, hv, and 
vv-curvature tensor components. For any Finsler space, the $hv$-curvature endomorphisms of the Chern-Rund connection are identically zero \cite{BCS}. However, for other linear connections in $\hat{\pi}^*TM$ the three curvature types could be different from zero.
\end{section}

\begin{section}{Riemannian average metrics from a Finsler space}
\subsection{Definition of the averaging procedure for Finsler metrics}
 The $n$-form $d^n y$ is defined in local coordinates by the expression
\begin{align}
d^n y=\sqrt{\det\, g(x,y)}\,\, \frac{\delta y^1}{F} \wedge \cdot \cdot \cdot \wedge \frac{\delta y^n}{F},
\label{volumeformontangentspace}
\end{align}
where $\det\, g(x,y)$ is the determinant of the fundamental tensor $(g)_{ij}=g_{ij}(x,y)$. Since
$d^n y$ is invariant by local diffeomorphisms on $TM$, it defines a section of $\Lambda^n N$.
 For each embedding $i_x:I_x \hookrightarrow N_x$ one consider the volume form on $I_x$ given by
\begin{align}
dvol_x:=i^*_x (d^n y\cdot\l ),
\label{volumeformontheindicatrix}
\end{align}
where $\l =\,\frac{y^i}{F}\,\frac{\partial}{\partial y^i}$ and $d^n y\cdot\l$ is the corresponding contraction.
Then the volume function $vol(I_x)$ is defined by the expression
\begin{align}
vol:M\to \mathbb{R}^+,\quad x\mapsto vol(I_x)=\int _{I_x} \, |\psi|^2(x,y)\,dvol_x,
\label{volumenoftheindicatrix}
\end{align}
where the weight factor $|\psi|^2:TM\setminus\{0\}\to \mathbb{R}^+$ is a homogenous of degree zero in $y$, positive, smooth function on the tangent bundle.
The average of a function $f\in \mathcal{F}(\mathcal{I})$ is defined by the expression
\begin{align*}
\langle f\rangle_{\psi} (x):=\frac{1}{vol(I_x)} \int _{I_x }\,|\psi|^2(x,y)\,f(x,y)\,dvol_x.
\end{align*}

Given the Finsler space $(M,F)$ let us consider the matrix with components
\begin{align}
h_{ij}(\psi,x):=\langle g_{ij}(x,y)\rangle_{\psi},\,i,j=1,...,n,
\label{averagedmetric}
\end{align}
for each point $x\in M$. Then we have,
\begin{proposicion}
Let $(M,F)$ be a Finsler space. Then $\{h_{ij}(\psi, x)\}^{n}_{i,j=1}$
are the components of a Riemannian metric
\begin{align}
h_\psi(x)=h_{ij}(\psi,x)\, dx^i\otimes dx^j,\quad i,j=1,...,n.
\label{averagedmetric2}
\end{align}
\label{propositionaveragemetric}
\end{proposicion}
\begin{proof} The average \eqref{averagedmetric} of a positive defined, real and symmetric $n\times n$ 
matrix is also a positive, real, symmetric matrix. To show this, first we note that
\begin{align*}
g_{x}(y)(\tilde{y},\tilde{y})=g_{ij}(x,y)\tilde{y}^i\tilde{y}^j\geq\,0,\,\,y,\tilde{y}\in\,T_xM.
\end{align*}
This is because  $g_{ij}(x,y)$ for each fixed $y\in \,T_x M$ is a positive defined scalar product in $T_x M$. The symmetric property follows in analogous way. The tensorial character of the expression \eqref{averagedmetric2} is also direct.
\end{proof}

\begin{definicion}
Two Finsler spaces $(M,F_1)$, $(M,F_2)$ are said to be $h$-equivalent if the corresponding average metrics $h_1$ and $h_2$ defined by \eqref{averagedmetric2} are the same.
\label{hequivalence}
\end{definicion}
 $h$-equivalence is an equivalence relation in the set of Finsler metrics $M_F$ defined over M. The equivalence class of $g$ is denoted by [g]. The coset space is defined by $M_F/\sim$. We call this equivalence relation $F$-equivalence.

From now we choose the function $|\psi|^2=1$, if anything else is not stated.  In this case the averaging is called {\it isotropic}.

\subsection{A coordinate-free formula for $h$}
Let us consider the fiber metric
\begin{align*}
\bar{g}=g_{ij}(x,y)\hat{\pi}^* dx^i \otimes \hat{\pi}^* dx^j,\quad   \xi \in \pi^{-1}_2 (u), \,u\in {\hat{\pi}^{-1}(x)}.
\end{align*}
\begin{proposicion}
Let $(M,F)$ be a Finsler space. Then the relation
\begin{align}
h(x)(X,X)= \,\frac{1}{vol(I_x)}\,\Big(\int _{I_x }  dvol_x\,
\,\bar{g} \,(\hat{\pi}^*|_u(X),\hat{\pi}^*|_u(X))\Big)
\end{align}
\label{proposiciononcoordenatefreemetric}
holds good for each vector field $X\in \,\Gamma\, TM$.
\end{proposicion}
\begin{proof} For each  point $x\in M$, $h(x)=\langle g_{ij}\rangle (x)\, d x^i |_x \otimes d x^j 
|_x $, $u\in I_x$ and vector field $X\in \,\Gamma\, TM$ one has
\begin{align*}
h(x) (X,X) & = \,\frac{1}{vol(I_x)}\,\Big(\int _{I_x }  dvol_x\,
\,{g}_{ij}(u)\Big)\,X^i\,X^j\\
& =\,\frac{1}{vol(I_x)}\,\Big(\int _{I_x }  dvol_x\,
\,\bar{g}_{ij}(u)\Big)\,X^i\,X^j\\
& =\,\frac{1}{vol(I_x)}\,\Big(\int _{I_x }  dvol_x\,
\,\bar{g} \,(\hat{\pi}^*|_u(X),\hat{\pi}^*|_u(X))\Big).
\end{align*}
\end{proof}
\end{section}

\begin{section}{General form of the average operation and applications}
\subsection{Average of a family of tensor automorphisms}
For each tensor $S_z \in {
T}^{(p,q)} _zM$ with $\, v\in
\hat{\pi}^{-1}(z)$ and $z\in { U}\subset M$, the following isomorphism is defined:
\begin{align*}
\hat{\pi}^*|_v :\hat{\pi}^{-1}(x)\to \pi^{-1}_1(v)  ,\quad S _z\mapsto  \hat{\pi}^*
_v S_z .
\end{align*}
Consider a family of fiber preserving automorphisms
\begin{align*}
A:=\big\{A_w:\hat\pi ^*_w T^{(p,q)}M \to \hat{\pi}^*_w T^{(p,q)}M,\, w\in\,I_x,\,x\in M\big\}.
\end{align*}
Then one can define the vector valued integral operation,
\begin{align}
\Big(\int _{I_x} \pi_2|_u  A_u  \hat{\pi}^* _u  \Big)\cdot\,S_x := \int _{I_x}\,
 \pi_2 ( \,A_u\,  \hat{\pi}^* _u \,S_x)\,dvol_x ,\quad S_x,
\in T^{(p,q)}_x M,
\label{fiberintegration}
\end{align}
where $u=(x,y)$ is a point on the fiber $I_x \hookrightarrow N_x$.
This operation is a fiber integration on a sub-manifold $I_x$ of $T_xM$ and with values of $T^{(p,q)}_xM$, for each $x\in\,M$. The tensor $S_x$ is pulled-back $\{\hat{\pi}^*_u S_x,\,u\in N_x\subset \,N\}$ such that the following diagram commutes,
\begin{displaymath}
\xymatrix{\hat{\pi}^*_uS_x\, \ar[d]_{\hat\pi_1} \ar[r]^{\hat\pi_2} &
S_x \ar[d]^{\pi}\\
I_x \ar[r]^{\hat{\pi}} & x}
\end{displaymath}
for each $x\in M,\,u\in\,\hat{\pi}^{-1}(x)$. The chain of compositions defining the integral operation \eqref{fiberintegration} is the following. Given $x\in M,\,u\in\hat{\pi}^{-1}(x)$ and $S_x\in\, T^{(p,q)}_XM$,
\begin{align}
x\mapsto \,S_x\mapsto \{\hat{\pi}^*_u S_x\}\mapsto \{A_u(\hat{\pi}^*_u S_x)\}\mapsto \{\pi_2(A_u(\hat{\pi}^*_u S_x))\}\mapsto \int_{I_x} \pi_2(A_u(\hat{\pi}^*_u S_x))\,dvol_x.
\label{operationofaveraging}
\end{align}
\begin{definicion}
 The average operator of the family of automorphism $A$ is the automorphism
\begin{align*}
\langle A\rangle_x  : {T }^{(p,q)}_x & M \to  T^{(p,q)}_x M\\
& S_x \mapsto \frac{1}{vol(I_x)}\Big(\int _{I_x}
 \pi
_2 |_u  A_u  \hat{\pi}^* _u \Big)\cdot\,S_x ,\,\quad u\in {\pi^{-1}_1(x)},\,\forall \, S_x
\in T^{(p,q)}_x M.
\end{align*}
\label{averagingoperation}
\end{definicion}
\begin{proposicion}
The average operator associated with a family of a linear of operators $A$ is a geometric operator on $M$.
 \end{proposicion}
  \begin{proof} Let us consider first sections of the bundle $\hat{\pi}^*TM$ and two local basis $\{e_i,\,i=1,...,n\}$ and $\{\tilde{e}_i,\,i=1,...,n\}$ of $T_xM$. Then we have
\begin{align*}
\langle A\rangle (S(x)) & =\frac{1}{vol(I_x)}\Big(\int _{I_x}
 \pi
_2 |_u  A_u  \hat{\pi}^* _u \Big)\cdot\,S(x)\\
& = \frac{1}{vol(I_x)}\Big(\int _{I_x}
 \pi
_2 |_u  A_u  \hat{\pi}^* _u \Big)\cdot\,{S}^i(x) e_i(x)\\
& = \frac{1}{vol(I_x)}\Big(\int _{I_x}
 \pi
_2 |_u  A_u  \hat{\pi}^* _u\,\tilde{S}^k(x)\,\tilde{e}_k(x)\,dvol_x\Big)\\
& = \frac{1}{vol(I_x)}\Big(\int _{I_x}
 \pi
_2 |_u  A_u  \hat{\pi}^* _u \Big)\cdot\,\tilde{S}^k(x)\tilde{e}_k(x).
\end{align*}
 One can extend this calculation to families of homomorphisms of $\hat{\pi}^*T^{(p,q)}_x M$.
\end{proof}
The averaging operation can be extended to families of local operators acting on sections $\hat{\pi}^*T^{(p,q)}M$ by applying the construction \eqref{operationofaveraging} in the {\it definition}  \ref{averagingoperation} to the evaluation of sections $\hat{\pi}^*_v(S):=\,(\hat{\pi}^* S)(v)$, where $S\in\,\Gamma\, T^{(p,q)}M$ and $v\in U_u$ is a point in the open set $U_u$ containing $u$.
\subsection{The average of linear connections on $\hat{\pi}^*TM$}
We can average linear connections on $\hat{\pi}^*TM$ following the above general method.
\begin{teorema}
 Let $\nabla$ be a linear connection of the vector bundle $\hat{\pi}^*TM\to N$. Then there
is defined an affine connection of $M$ determined by the
covariant derivative of each section $Y\in\,\Gamma\, TM$ along each  direction $X\in T_x M$,
\begin{align}
\langle\nabla\rangle_X  Y:= \langle \pi _2|_u
{\nabla}_{\iota_u({X})} \hat{\pi}^* _v Y\,\rangle ,\,v\in TU_x\setminus 0,
\label{averagedconnectionformula}
\end{align}
for each $X\in T_x M$ and $Y\in\,\Gamma\, TM$, where ${U}_x$ is an open neighborhood of $x\in M$.
\label{averagedconnection}
\end{teorema}
\begin{proof}
We check that the properties for a linear covariant
derivative hold for $\langle \nabla\rangle$.
\begin{enumerate}

\item
Using the linearity of the original covariant derivative and the linearity of the averaging operation,
\begin{align*}
\langle\nabla\rangle_{X }(Y_1 +Y_2 ) &=\langle\pi _2|_u {\nabla}_{{\iota}_u (X)} \hat{\pi}^* _v 
(Y_1 +Y_2 )\rangle & =\langle\pi _2 |_u  {\nabla}_{{\iota}_u (X)} \hat{\pi}^* _v Y_1
\rangle +\langle\pi _2 |_u   {\nabla}_{{\iota}_u (X)}  \hat{\pi}^* _v Y_2 \rangle\\
& =\langle\nabla\rangle_{X} Y_1 +\langle\nabla\rangle_{X} Y_2,
\end{align*}
$\forall \,\,  Y_1 , Y_2$.
For the second condition of linearity we have
\begin{align*}
\langle\nabla\rangle_{X } ({\lambda}Y)=\langle\pi _2 |_u {\nabla}_{{\iota}_u
(X)}
\hat{\pi}^*_v ({\lambda}Y)\rangle_v ={\lambda}\langle\pi _2 |_u  {\nabla}_{{\iota}_u (X)} 
\hat{\pi}^*_v \rangle=\,\lambda\langle\nabla\rangle_{X } Y
\end{align*}
 $\forall\,\, Y\in { \Gamma TM}, \lambda\in\,\mathbb{R}, \, X\in T_x M.$
\item $\langle\nabla\rangle_{X } Y$ is a ${\mathcal{F}}$-linear
respect $X$, that is,
\begin{align}
\langle\nabla\rangle_{X_1 +X_2 } Y=\langle\nabla\rangle_{X_1 }
Y+\langle\nabla\rangle_{X_2 } Y,\quad\quad \langle\nabla\rangle_{fX} (Y)=f\,\langle\nabla\rangle_{X } Y,
\end{align}
holds good for each $Y\in\,\Gamma T M,v\in {\pi}^{-1}(z),\, X,X_1,X_2 \in T_x M$ and $f\in{\bf \mathcal{F}}M.$
The first condition is proved by the following short calculation,
calculation,
\begin{align*}
\langle\nabla\rangle_{X_1 +X_2 }Y & =\langle \pi _2 |_u ({\nabla}_{{\iota}_u
(X_1 +X_2 )})\hat{\pi}^*_v Y\rangle_v =\langle\pi _2 |_u {\nabla}_{{\iota}_u (X_1 )}\hat{\pi}^*_v Y \rangle_v +\langle\pi _2
|_u {\nabla}_{{\iota _u}(X_2 )}\hat{\pi}^*_v Y \rangle\\
& =\langle\nabla\rangle_{X_1 } Y+\langle\nabla\rangle_{X_2 } Y.
\end{align*}
For the second condition the proof is similar.

\item The Leibnitz rule holds:
\begin{align}
\langle\nabla\rangle_{X }(\varphi Y)=\,(d\varphi(X))Y+\varphi\langle\nabla\rangle_{X } Y,\quad \forall \,  Y\in \,\Gamma TM,\quad \varphi\in \mathcal{F}(M),\quad X \in T_x M,
\label{Leibnitz}
\end{align}
 where $d\varphi(X)$ is the action of the $1$-form $d\varphi\in { \Lambda}^1  M$ evaluated at $x\in M$ on the tangent vector $X\in 
T_x M$. In order to prove \eqref{Leibnitz} we use the following property,
\begin{displaymath}
\hat{\pi}^*_v (\varphi Y)=\,(\hat{\pi}^*_v \varphi)\,(\hat{\pi}^*_v Y),\quad \forall \,\, Y\in {
TM},\, \varphi\in {\mathcal{F}}(M).
\end{displaymath}
Then one obtains the following expressions,
\begin{align*}
\langle{\nabla}\rangle_{X} (\varphi Y) & =\langle \pi _2 |_u {\nabla}_{{\iota}_u
(X)}\hat{\pi}^*_v (\varphi Y)\rangle =\langle \pi _2 |_u {\nabla}_{{\iota}_u (X)}
(\hat{\pi}^*_v \varphi)\hat{\pi}^*_v Y\rangle\\
& =\,\langle \pi _2 |_u (\nabla _{{\iota}_u (X)}( \hat{\pi}^* _v \varphi))\hat{\pi}^*_v (Y)\rangle+\langle\pi _2 
|_u (\hat{\pi}^* _u f){\nabla}_{{\iota}_u (X)}\hat{\pi}^*_v (Y)\rangle\\
& =\,\langle\pi _2 |_u  ({{\iota}_u (X)}\cdot(\hat{\pi}^*_v \varphi)) \hat{\pi}^*_v (Y)\rangle + \varphi\, \langle\pi _2 |_u {\nabla}_{{\iota}_u (X)} \hat{\pi}^*_v (Y)\rangle \\
& =\,\langle (X\cdot \varphi)\pi _2 |_u \hat{\pi}^*_u (Y)\rangle +\varphi\, \langle\pi _2 |_v
{\nabla}_{{\iota}_u (X)}\hat{\pi}^*_v (Y)\rangle.
\end{align*}
For the first term we perform the following simplification,
\begin{align*}
\langle(X\cdot \varphi)\pi _2 |_u \hat{\pi}^*_u (Y)\rangle  & =\,(X \cdot \varphi)\langle \pi _2 |_u  \hat{\pi}^*_u 
(Y)\rangle =\,(X \cdot\varphi)(\langle  \pi _2 |_u \hat{\pi}^*_u \rangle )Y  \\
& =\,(X\cdot \varphi)Y =\,(d\varphi(X))Y.
\end{align*}
Finally we obtain that
\begin{align*}
\langle{\nabla}\rangle_{X} (\varphi Y) =\,(\langle \nabla\rangle _{X_x}\varphi)Y+\varphi
\langle{\nabla}\rangle_{X}Y =\,(d\varphi(X))Y+\varphi\langle\nabla\rangle_{X }Y.
\end{align*}
\end{enumerate}
\end{proof}
We denote the affine connection associated with the above covariant derivative by 
$\langle{{\nabla}}\rangle$. Then for each section $Y\in \Gamma \,TM$, $\langle \nabla \rangle Y\in \Gamma \,T^{(1,1)}
M$ is given by
\begin{equation}
\langle \nabla \rangle(X,Y):=\langle \nabla \rangle_X Y,\quad X\in T_xM.
\end{equation}
The average covariant derivative is extended to $1$-forms by the requirement that it commutes with contractions.
 Thus for each $\alpha\in\, \Lambda^1 M$ and $X\in \, \Gamma TM$
\begin{align*}
\langle \nabla \rangle_X [\alpha (Z)]=\,(\langle\nabla \rangle_X \alpha)\cdot Z
+\, \alpha \cdot(\langle \nabla \rangle_X Z)
\end{align*}
holds by assumption.
Then the extension of the covariant derivative $\langle \nabla \rangle_X$ to sections of the tensor bundle $T^{(p,q)}M\to M$ is defined by the rule
\begin{align*}
\langle \nabla \rangle_X
\,K(X_1,...,X_s,\alpha^1,...,\alpha^r) &=\langle \nabla \rangle_X\,
K(X_1,...,X_s,\alpha^1,...,\alpha^r)\\
& -\sum^s_{i=1}K(X_1,...,\langle \nabla \rangle_X
X_i,...,X_s,\alpha^1,...,\alpha^r) \\
& +\sum^r_{j=1}
K(X_1,...,X_s,\alpha^1,...,\langle \nabla \rangle_X \alpha^j,...,\alpha^r).
\end{align*}
\subsection{General properties of the averaged connection}
\begin{proposicion}
Let $(M,F)$ be a Finsler space and $\nabla$ a linear connection on $\hat{\pi}^*TM$. Then it holds
\begin{align}
T_{\langle \nabla \rangle}=\,\langle T_\nabla\rangle.
\label{averagetorsion}
\end{align}
\end{proposicion}
\begin{proof}
We can calculate the torsion of the connection $\langle \nabla \rangle$,
\begin{align*}
T_{\langle \nabla \rangle}(X,Y) & =\langle\pi _2 |_u \nabla _{\iota _u (X)} \hat{\pi}^* _v \rangle Y -\langle\pi 
_2 |_u  \nabla _{\iota _u (Y)} \hat{\pi}^* _v \rangle X -[X,Y]\\
& =\langle\pi _2 |_u \nabla _{\iota _u (X)}  \hat{\pi}^* _v \rangle Y -\langle\pi _2 |_u 
\nabla _{\iota _u (Y)} \hat{\pi}^* _v \rangle X -\langle\pi _2 |_u  \hat{\pi}^* _u [X,Y]\rangle\\
& =\langle\pi _2 |_u  \big(\nabla _{\iota _u (X)}\hat{\pi}^* Y -\nabla _{\iota _u (Y)}\hat{\pi}^*X 
-\hat{\pi}^*|_u[X,Y] \big) \rangle\\
& =\langle T_\nabla (X,Y)\rangle.
\end{align*}
\end{proof}
\begin{corolario}
Let $(M,F)$ be a Finsler space with average Chern-Rund connection $\langle\,^{ch}\nabla \rangle$. Then the torsion $T_{\langle \,^{ch}\nabla \rangle}$  is zero.
\label{averageofthetorsionforchern}
\end{corolario}
\begin{definicion}
Two Finsler spaces $(M,F_1)$ and $(M,F_2)$ are $\Gamma$-related iff the corresponding average connections are the same.
\end{definicion}
This is an equivalence relation ($\Gamma$-equivalence). $\Gamma$-equivalence classes are denoted as $[g]_\Gamma$. Let us remark that the averaging procedure depends on the particular Finsler spacer $(M,F)$ that one starts and that this equivalence relation is between different types of averaging procedures. Also note that in general the $\Gamma$-equivalence relation is different than the $F$-equivalence relation.

\subsection{Average connection of a Berwald space}
In order to avoid cluttering in the notation, the Chern-Rund connection will be denoted simply by $\nabla$ instead of the most specific notation $^{ch}\nabla$ used until now.
\begin{definicion}
A Berwald space is a Finsler space such that its Chern-Rund connection also defines  an affine connection on $M$.
 \end{definicion}
For a Berwald space the connection coefficients $^{ch}\Gamma^i_{jk}(x,y)$ depend on $x\in\,M$ only. Thus, we have the following,
\begin{teorema}
For a Berwald space $(M,F)$
\begin{itemize}
 \item The average of the Chern-Rund connection coincides with the Chern-Rund connection in the sense that
\begin{align}
\left(\hat\pi^*\langle \nabla\rangle_{X}\,S\right)=\nabla_{\iota_u(X)}\hat{\pi}^*S.
\label{averageconnection=connectioninBerwald}
\end{align}
\item The average of the Chern-Rund connection coincides with the Levi-Civita connection of the average metric,
\begin{align}
\langle \nabla\rangle\,=\,^h\nabla.
\label{averageconnection=LeviCivitaconnectioninBerwald}
\end{align}
\end{itemize}
\label{averageconnectionberwaldspaces}
\end{teorema}
\begin{proof} The relation \eqref{averageconnection=connectioninBerwald} is direct from the definition of average connection. A detailed proof can be found in \cite{RF}. For the proof of  relation \eqref{averageconnection=LeviCivitaconnectioninBerwald}, let us  first evaluate the covariant derivative of the metric $h$ for $\langle \nabla\rangle$ ,
\begin{align*}
\langle \nabla\rangle_{X=\frac{\partial}{\partial x^i}}\,h &=\,\Big(\frac{\partial h_{jk}}{\partial x^i}- \,\langle \nabla\rangle^l\,_{ik}h_{jl}-\,\langle \nabla\rangle^l\,_{ij}h_{lk}\Big)dx^j\otimes dx^k\\
& =\,\Big(\frac{\partial h_{jk}}{\partial x^i}- \,^{ch}\Gamma^l\,_{ik}h_{jl}-\,^{ch}\Gamma^l\,_{ij}h_{lk}\Big)dx^j\otimes dx^k\\
& = \frac{1}{vol(I_x)}\,\Big(\int_{I_x}\,\big(\frac{\partial g_{jk}}{\partial x^i}-\,^{ch}\Gamma^l\,_{ik}g_{jl}-\,^{ch}\Gamma^l\,_{ij}g_{lk}\big)\,dvol_x\Big)dx^j\otimes dx^k\\
& =\frac{1}{vol(I_x)}\,\Big(\int_{I_x}\big(\nabla_{\frac{\partial}{\partial \tilde{x}^i}}\,g)\big)\,_{jk}\,dvol_x\Big)dx^j\otimes dx^k.
\end{align*}
Since the horizontal metric compatibility of the Chern-Rund connection ({\it equation} \eqref{horizontalmetricondition}), the integrand is zero. Therefore,
\begin{align*}
\langle \nabla\rangle_{X=\frac{\partial}{\partial x^i}}\,h=0,\,\quad i=1,...,n.
\end{align*}
Moreover, by {\it Corollary} \ref{averageofthetorsionforchern} $\langle \nabla\rangle$ is torsion free. Therefore, $\langle \nabla\rangle$  must be the Levi-Civita connection of $h$.
\end{proof}

The geodesic deviation equation for a Finsler space is formally the same than for a Riemannian space. In particular, the Jacobi equation for the Chern-Rund connection \cite{BCS} is the linear differential equation along the geodesic $X:I\to M$
\begin{align}
\nabla_X\nabla_X J+R(X,J)X=0.
\label{Jacobi equation for a Finsler metric}
\end{align}
It is interesting that if the space $(M,F)$ is of Berwald type, then the equation \eqref{Jacobi equation for a Finsler metric} is formally the same than the equation for the geodesics of the average connection,
\begin{align*}
\langle\, \nabla\rangle_{\tilde {X}}\langle\, \nabla\rangle_{\tilde{X}} \tilde{J}+R^{\langle \nabla \rangle}(\tilde{X},\tilde{J})\tilde{X}=0,
\end{align*}
or
\begin{align*}
\nabla_{\tilde{X}}\nabla_{\tilde{X}} \tilde{J}+R(\tilde{X},\tilde{J})\tilde{X}=0,
\end{align*}
where $\tilde{X}:\tilde{I}\to M$ is a geodesic of the averaged connection $\langle \nabla\rangle$ and $R^{\langle \nabla \rangle}$ is the curvature endomorphism.
The Jacobi fields $J:I\to X(I)$ and $\tilde{J}:\tilde{I}\to \tilde{X}(\tilde{I})$ are the same (at least for small values of the  time values of the time parameter,  the geodesics $X$ and $\tilde{X}$ as un-parameterized geodesics), except that they are parameterized by different parameters. Then  the condition
\begin{align*}
R^{\langle\nabla \rangle}(\tilde{X},\tilde{J})=R(X,J)
\end{align*}
must hold for Berwald spaces.
Since for a Berwald space $\langle \nabla\rangle=\,^h\nabla$, the following result holds:
\begin{teorema}
If $(M,F)$ is a Berwald space, then
\begin{align}
^hR(X,Y)=R(X,Y),\quad X,Y\in \,\Gamma TM.
\label{relation R endomorphisms}
\end{align}
\label{teorema on the relation between R endomorphism}
\end{teorema}
\begin{corolario}
If $(M,F)$ is a Berwald space, then $^hR^i\,_{jkl} =\,\langle\,^gR^i\,_{jkl}\rangle$.
\end{corolario}
\begin{proof}
This is a direct consequence of the above calculation and {\it Theorem} \ref{teorema on the relation between R endomorphism}.
\end{proof}

A simple calculation shows that in general the average of the curvature endomorphism does not coincide with the curvature endomorphism of the average metric,
\begin{align*}
^hR^i\,_{jkl} & =\,h^{im}\,^hR_{mjkl}  =\,h^{im}\,\langle\,^gR_{mjkl}\rangle
=\,\langle\,\,h^{im}\,^gR_{mjkl}\rangle\\
& =\langle\,h^{im}g_{ms}g^{sa}\,^gR_{ajkl}\rangle=\,\langle\,\vartheta ^i \,_s\,^gR^s\,_{jkl}\rangle,
\end{align*}
where the tensor $\vartheta^i\,_s:=h^{im}g_{ms} \neq \delta^i_s$ measures the departure of the fundamental tensor $g$ of being Riemannian.
Similarly, in the general case there is no direct relation between the average of the flag curvature of $F$ and the sectional curvatures of $h$, even in the case of constant sign flag curvature spaces.

Let $(M,F)$ be a Berwald space. Then  the Riemann tensor of $g$ is a $(0,4)$-tensor along the map $\hat{\pi}:N\to M$ whose components are given in normal coordinates of $g$ by the expression\footnote{For a Berwald space, normal coordinate system exists and are $\mathcal{C}^2$, see for instance \cite{BCS}.}
\begin{align}
^gR_{ijkl}:=g_{il,jk}\,- g_{ik,jl}\,+ g_{jk,il}\,- g_{jl,ik},
\label{riemanntensorofg}
\end{align}
where $g_{ij,kl}$ stands for $\frac{\partial^2 g_{ij}}{\partial x^k\partial x^l}$, etc...
\begin{proposicion}
Let $(M,F)$ be a Berwald space and consider the isometric average metric $h_{ij}$. Then
the following relation holds,
\begin{align}
^{h}R_{ijkl}(x)=\,\langle\,^gR_{ijkl}(x,y)\rangle.
\label{averageofh}
\end{align}
\label{propositiononthecurvatureoftheaverage}
\end{proposicion}
\begin{proof} In normal coordinates for $h$, the Riemann tensor $^h R_{ijkl}$  is linear on the second derivatives of the components of the curvature tensor of $h$. Then the tensor can be expressed as
\begin{align*}
^hR_{ijkl}=\,h_{il,jk}\,-h_{ik,jl}\,+h_{jk,il}\,-h_{jl,ik},
\end{align*}
where for instance $h_{il,jk}=\frac{\partial^2 h_{il}(x)}{\partial x^j\partial x^k}$, etc...
From the definition of $h$ it follows that
\begin{align*}
^hR_{ijkl}=\,\langle\,g_{il}\rangle,_{jk}\,-\langle g_{ik}\rangle,_{jl}\,+\langle g_{jk}\rangle,_{il}\,-\langle h_{jl}\rangle,_{ik}
\end{align*}
holds in the normal coordinate chard of $h$.
Since the weight factor is $|\psi|^2=\,{1}$ and the volume function $vol(I_x)$ is constant for a Berwald space \cite{BaoShen}, the partial derivatives can be introduced in the integrals,
\begin{align*}
^{\langle g\rangle}R_{ijkl}(x) & =\,^hR_{ijkl}=\,\langle\,g_{il}\rangle,_{jk}\,-\langle g_{ik}\rangle,_{jl}\,+\langle g_{jk}\rangle,_{il}\,-\langle h_{jl}\rangle,_{ik}\\
& =\,\langle\,g_{il,jk}\,- g_{ik,jl}\,+ g_{jk,il}\,- g_{jl,ik}\rangle=\,\langle \,^gR_{ijkl}(x,y)\rangle,
\end{align*}
where from the first to the second line in the above expression we have use the equality $^{ch}\Gamma^i\,_{jk}(x)=\,^h\Gamma^i\,_{jk}(x)=0$ in normal coordinates of the average metric $h$. In  this step is essential the Berwald condition.
This formulae has been proved in normal coordinates of $h$, but since it is an identity between tensor components, it holds in any coordinate system.
\end{proof}
As an example of direct application of {\it proposition} \ref{propositiononthecurvatureoftheaverage} is a generalization of Riemann's characterization of Euclidean space in terms of curvature. Let $h_o$ be the Euclidean metric in $\mathbb{R}^n$. Then we have
\begin{corolario}
Let $(M,F)$ be a Berwald space such that $^gR_{ijkl}=0$. Then $M\cong \mathbb{R}^n$ and the fundamental tensor $g=\,h_0+\delta g$, with $\langle\delta g\rangle =\,0$.
\end{corolario}
\begin{proof}
If $^gR_{ijkl}=0$, then the relation \eqref{averageofh} implies $^{h}R_{ijkl}(x)=0$ and the result follows from the classical result of Riemann \cite{Riemann}.
\end{proof}

\subsection*{Gauss-Bonnet theorem for Berwald surfaces}The construction of the  metric $h$ as an average of the fundamental tensor over the indicatrix opens the possibility to generalize results from Riemannian to Berwald geometry, using directly the Riemannian results. We consider here another example of this technique: a weak version of the Gauss-Bonnet theorem for arbitrary Berwald surfaces.  Let $(\mathcal{I},\pi_{\mathcal{I}},M)$ be the fibered manifold whose fibers are indicatrix over $M$ and note that for a Berwald space the function $x\mapsto vol(I_x)$ is constant. Then we have
\begin{teorema}
Let $(M, F)$ be a compact  Berwald surface with average metric $h$ and Gaussian curvature $^h K=\,-\,^hR_{1212}$ in some orthonormal basis of $h$. Then the following formula holds \cite{Dazord},
 \begin{align}
 \frac{1}{vol(I_x)} \,\int_{\mathcal{I}}\, ^gR_{1212}(x,y)\,dvol_x\wedge d\mu(x)=\,-2\pi\,\chi(M),
 \label{GaussBonnet2}
 \end{align}
 where $\chi(M)$ is the Euler's characteristic of $M$ and $d\mu$ is the Riemannian volume form of $h$ on $M$.
 \label{teoremagaussbonet}
\end{teorema}
\begin{proof} For the Riemannian metric $h$ one can apply the classical Gauss-Bonnet theorem for compact surfaces $M$. Thus the relation
\begin{align}
\int_M \, ^hR_{1212}(x)\,d\mu=\, -2\pi\,\chi(M)
\label{GaussBonnet1}
\end{align}
holds.
Fixed the integration measure as in {\it proposition} \ref{propositiononthecurvatureoftheaverage} an using an orthonormal frame associated to $h$, one obtains the relations
\begin{align*}
^hK(x)=-\,^hR_{1212}(x)=-\langle \,^g R_{1212} \rangle,
\end{align*}
 from which the relation \eqref{GaussBonnet2} follows.
\end{proof}

\section{The parallel transport of the average connection}
 The parallel transport of a linear connection ${\nabla}$ along a path $\gamma:[a,b]\to M$
  with $\gamma(a)=x$ and $\gamma(b)=z$ is defined as the linear homomorphism
 \begin{align*}
 p_{xz}(\gamma)  :T^{(p,q)}_x &M \rightarrow T^{(p,q)}_zM,\quad S_x\mapsto S_z
 \end{align*}
such that the section $S(t)$ along $\gamma$ is a solution of the linear differential equation
\begin{align}
{\nabla}_{\dot{\gamma}(0)}p(\gamma)(S)(t)=0,\quad p(\gamma)(S)(0)=S(0).
\label{equationparalleltransport}
\end{align}
A {\it polygonal approximation} $\bar{\gamma}$ of $\gamma$ is determined by a set of points
\begin{align*}
\{\gamma(0)=x,\,...,\gamma(t_i),\,...,\,\gamma(t_{A-1}),\gamma(t_A)=z,\,\gamma({t_i})\in  \, \gamma([a,b])\}
\end{align*}
 joined by geodesic segments $\bar{\gamma}_{k,k+1}$ of $F$, with initial and ending points $\gamma(t_{k})$ and $\gamma(t_{k+1})$ respectively. One can also consider the case when $t_k -t_{k-1}=\,\epsilon$. Then the  parallel transport operator along $\bar{\gamma}$ is given by the composition of  parallel transports
\begin{displaymath}
p(\bar\gamma):=\,\prod^A_{k=1}\,\circ\, p_{t_k,t_{k-1}},
\end{displaymath}
where the composition of elementary parallel transport $p_{t_k,t_{k-1}}$ is taken along the geodesic segment $\bar{\gamma}_{k,k+1}$  and is given by the endomorphism
 \begin{align*}
 p_{t_k,t_{k-1}}  :& T_{\gamma(t_{k-1})}M \to T_{\gamma(t_{k})}M\\
  & X^ie_i\mapsto \big(\delta^{j_k}_{i{_{k-1}}} X^{i_{k-1}} \,-\epsilon\,{\Gamma}^{j_k}_{i_{k-1}l_ {k-1}}\,X^{i_{k-1}}\dot{\gamma}^{l_{k-1}}\big)e_ {j_{k}}.
 \end{align*}
$\dot{\gamma}^{l_{k-1}}$ is the tangent vector at the point $\gamma(k-1)$.
Then the double limit $A\rightarrow +\infty$ and $\epsilon =t_k-\,t_{k-1}\rightarrow 0$ is taken in this parallel transport operation, under the constraint
\begin{align*}
\lim_{A\to +\infty ,\,\epsilon\to 0}\,A\epsilon\,=b-a.
\end{align*}
Let us define $\epsilon=\,\frac{b-a}{A}$.
Then the  parallel transport of $X\in \,T_{\gamma(0)}M$ along $\gamma$ between the point $x=\gamma(a)$ and $z=\gamma(b)$ is given by
\begin{align}
 (p_{xz}X)^j\,=\,\lim_{A\rightarrow +\infty,\epsilon\rightarrow 0}\,\Big( (\delta^{j_A}_{i{_{A-1}}}   \,-\epsilon\,\Gamma^{j_A}_{i_{A-1}l_
 {A-1}}\,\dot{\gamma}^{l_{A-1}})\Big) (p_{x\gamma(t_{A-1})}X)^{i{_{A-1}}},\,j=1,...,n,
 \label{paralleltransportformula}
\end{align}
with $(p_{xx}X)^{j_0}=\,X^{j_0}$, $\lim_{A\rightarrow +\infty}\gamma(t_{A-1})=\gamma(b)$ and $j_A=j$.

 Infinitesimally one has the finite difference expression
  \begin{align*}
 (p_{x{\gamma}(t+\epsilon)}X)^j-\,(p_{x{\gamma}(t)}X)^j & =\,\Big( (\delta^{j_A}_{i{_{A-1}}}   \,-\epsilon\,\Gamma^{j_A}\,_{i_{A-1}l_
 {A-1}}\,\dot{\gamma}^{l_{A-1}})\Big) (p_{x{\gamma}(t)}X)^{i{_{A-1}}}\\
 & \,-(p_{x\gamma(t)}X)^j\\
 & =\,-\,\epsilon\,(p_{x{\gamma}(t)}\Gamma^j_{ik}(\gamma(t))\,\dot{\gamma}^i(t)(p_{x\gamma(t)}X)^k.
\end{align*}
For smooth vector fields and connections, one can take the limit
\begin{align*}
 \lim_{\epsilon\to 0}\,\frac{(p_{x\gamma(t+\epsilon)}X)^j-\,(p_{x\gamma(t)}X)^j}{\epsilon}
  =\,-(p_{x\gamma(t)}\Gamma^j_{ik}(\gamma(t))\,\dot{\gamma}^i(t)(p_{x\gamma(t)}X)^k,
\end{align*}
showing that the expression \eqref{paralleltransportformula} is the solution of the parallel transport equation \eqref{equationparalleltransport}.

The expression \eqref{paralleltransportformula} applies to the parallel transport of any linear connection. In particular, it can be applied to the average connection,
\begin{proposicion}
Let $\langle \Gamma^i_{jk}\rangle $ be the connection coefficients of the average connection $\langle \,\nabla\rangle $. Then the parallel transport operation is
\begin{align}
(p_{xz}X)^j\,=\,\lim_{A\rightarrow \infty,\,\epsilon\rightarrow 0}\,\prod^A_{k=1}\Big( (\delta^{j_k}\,_{i{_{k-1}}} \,-\epsilon\,\langle \Gamma^{j_k}_{i_{k-1}l_
{k}}\rangle ({\gamma(t_{k})}) \,\dot{\gamma}^{l_{k}})\Big) X^{i_0},\,\quad j=1,...,n.
\label{averageparalleltransportformula}
\end{align}
\end{proposicion}
If we explicitly insert the weight factor $|\psi|^2:N\to {\mathbb{R}}^+$ in each integration, one obtains the expression for the parallel transport
\begin{align}
& (p_{xz}X)^j  \,=\,\lim_{A\rightarrow +\infty,\epsilon\rightarrow 0}\,\Big(\prod^A_{k=1}
\frac{1}{vol(I_{\gamma(t_{k-1})})}\,\Big(\int_{I_{\gamma(t_{k-1})}}\,dvol_x(\gamma(t_{k-1}))\, \\
& \nonumber f(\gamma(t_{k}),y_{\gamma(t_{k})})\,(\delta^{j_k}\,_{i{_{k-1}}}\,-\epsilon\,\Gamma^{j_k}\,_{i_{k-1}l_
{k-1}}(\gamma(t_{k}),y_{\gamma(t_{k})})\,\dot{\gamma}^{l_{k-1}}(\gamma(t_{k-1}))\Big)\Big) X^{i_0}.
\end{align}

 In the general case the average of the parallel transport operation does not coincide with  the parallel transport of the average connection:
the average of the parallel transport of $\nabla$ involves only one fiber integration, while the parallel transport of the average connection involves
a formal infinite number of integral operations along each fiber $\hat{\pi}^{-1}(\gamma(t)),\,t\in\,[a,b]$.

\subsection{Curvature of the average connection}
Let us consider  the curvature endomorphisms for the average connection $\langle\nabla\rangle$,
\begin{displaymath}
R^{\langle\nabla\rangle}(X_1,X_2)Z=\,\big(\langle\nabla\rangle_{X_1}\langle\nabla\rangle_{ X_2}
\,-\langle\nabla\rangle_{X_2}\langle\nabla\rangle_{X_1}\,-\langle\nabla\rangle_{ [X_1,X_2]}\big)Z.
\end{displaymath}
Developing this expression in terms of the original connection $\nabla$ one obtains
\begin{align*}
 R^{\langle\nabla\rangle}_x(X_1,X_2)(Z) & =\,\frac{1}{vol^2(\Sigma_{x})}\,\int_{I_{x}}\,\int_{I_x}\, dvol_{x}(v)dvol_{x}(u)\,\pi_2(v)\\
 &\Big(\,\nabla_{\iota_v(X_1)}\hat{\pi}^*_v\,\pi_2(u)\,\nabla_{\iota_u(X_2)}
 \hat{\pi}^*_u\,Z-\,\nabla_{\iota_v(X_2)}\hat{\pi}^*_v\,\pi_2(u)\,\nabla_{\iota_u(X_1)}\hat{\pi}^*_u\,Z\\
& -\,\nabla_{\iota_u([X_1,X_2])}\,\hat{\pi}^*_u \,Z\,\Big).
\end{align*}

It is interesting that the  curvature $R^{\langle \nabla \rangle}$ is not equal to the {\it average curvature} of the linear connection $\nabla$. For instance, the averaged $hh$-curvature is
\begin{align*}
\langle R^{\nabla}(\,\iota_u({X}_1),\iota_u(X_2))\rangle Z\, & :=\,\frac{1}{vol(I_x)}\,\int_{I_x}\, dvol_{x} \,\pi_2(u)\,\nabla_{\iota_u(X_1)}\nabla_{\iota_u(X_2)}\hat{\pi}^*_u\,Z\\
& -\frac{1}{vol^(I_x)}\,\int_{I_x}\,dvol_{x} \,\pi_2(u)\nabla_{\iota_u(X_2)}\,\nabla_{\iota_u(X_1)}\hat{\pi}^*_u\,Z\\
& -\frac{1}{vol(I_x)}\,\int_{I_x}\,dvol_{x} \,d u\, \pi_2(u)\,\nabla_{\iota_u([X_1,X_2])}\,\hat{\pi}^*_u\,Z\\
& =\frac{1}{vol(I_x)}\,\int_{I_x}\, dvol_{x}
\,\pi_2(u)\,\Big(\nabla_{\iota_u(X_1)}\nabla_{\iota_u(X_2)}\\
& -\nabla_{\iota_u(X_2)}\,\nabla_{\iota_u(X_1)}\,-\nabla_{\iota_u([X_1,X_2])}\,\Big)\hat{\pi}^*_u\,Z.
\end{align*}
Therefore, given a linear connection on $\pi^*TM$ $\nabla$, there are two notions of {\it average curvature endomorphisms}, $R^{\langle\nabla\rangle}_x(X_1,X_2)$ and $\langle R^{\nabla}(\,\iota_u({X}_1),\iota_u(X_2))\rangle$. In the general case the tensors $R^{\langle\nabla\rangle}_x(X_1,X_2)$ and $\langle R^{\nabla}(\,\iota_u({X}_1),\iota_u(X_2))\rangle$ do not coincide because the covariant derivative $\nabla_{\iota_u(Y)}$ depends on $u\in N$ for a general Finsler space.

\subsection{Holonomy of a Berwald space}
For a Berwald space, the Chern-Rund connection lives on the manifold $M$ and the curvatures of the average connection coincide with the average of the curvature of the original linear connection on $\pi^*TM$. Therefore, as an application of the Ambrose-Singer theorem on holonomy \cite{Ambrose Singer, Kobayashi Nomizu} and {\it theorem} \ref{teorema on the relation between R endomorphism} it holds the following result. Let us consider the holonomy group $Hol(\nabla)$ of the Chern-Rund connection $\nabla$. Then
\begin{teorema}
Let $(M,F)$ be a Berwald space. Then the holonomy group $Hol(\nabla)$ is Riemannian.
\end{teorema}
\begin{proof}
If the space $(M,F)$ is a Berwald space, then the $hv$-curvature endomorphisms are identically zero, $P=0$. Then the result follows from a direct application of the Ambrose-Singer theorem on holonomy and the relation \eqref{relation R endomorphisms}.
\end{proof}
An interesting consequence is that, if the space $(M,F)$ is Berwald, then the average holonomy group is not only an affine holonomy group, but it is indeed metrizable. This is an extension of a theorem from Z. Szab\'o on metrizability of compact holonomies of Berwald spaces \cite{Szabo1980},
\begin{teorema}
Let $(M,F)$ be a Berwald space. Then the holonomy of the Chern-Rund connection is metrizable.
\label{generalizacion del teorema de szabo}
\end{teorema}
\section{Isometries of the average metric}
\begin{definicion}
 Given two Finsler spaces $(M_1, F_1)$ and $(M_2,F_2)$,
a {\it base manifold Finsler isometry} (or simply a Finsler isometry) is a diffeomorphism $\Phi:M_1\to M_2$ such that preserves the Finsler function,
\begin{align}
F_2(\Phi (x),d\Phi(y))=F_1 (x,y).
\end{align}
\label{definicionfinslerisometry}
\end{definicion}
As a direct consequence of this definition and in the case when $M_1=\,M_2$, the components of the fundamental tensor transform locally under an isometry as
\begin{align}
{ (g_1)}_{ij}(\tilde{x}(x),\tilde{y}(x,y))=\frac {\partial x^l }{\partial
\tilde{x}^i}\frac{\partial x^k }{\partial \tilde{x}^j}{(g_2)}_{lk}(x,y).
\label{local map transformation under isometry}
\end{align}

It turns out that \cite{GallegoTorromePiccione}  
\begin{proposicion}
The group of isometries of $(M,F)$ is contained as a closed subgroup of the isometries of $(M,h)$ (in the compact-open topology).
\label{iosmetry of h and isometry of F}
\end{proposicion}
Then it follows that the group of isometries\footnote{Our result is a particular application of the averaged method. Other example is found in \cite{Matveev Troyanov 2017}, where a different average method was applied \cite{Matveev Troyanov 2012}.} of $(M,F)$ is a Lie group \cite{DengHou02} and it is a subgroup of the group of isometries of the average metric $(M,h)$.

\begin{definicion}
A Finsler space $(M,F)$ is symmetric if for each point $x\in M$ there is a Finsler isometry $\varphi_x:M\to M$ such that
\begin{itemize}
\item $\varphi_x(x)=x,$

\item $(d\varphi_x)|_x=\,-Id|_{T_xM}$.
\end{itemize}
\end{definicion}
The two conditions of the isometry $\varphi$ are the same for $h$ than for $F$. Then  by direct application of {\it proposition} \eqref{iosmetry of h and isometry of F} we have
\begin{teorema}
If $(M,F)$ is a Finsler symmetric space, then $(M,h)$ is a Riemannian symmetric space.
\end{teorema}
This result holds for locally symmetric spaces and globally symmetric spaces.
In particular, we have
\begin{corolario}
If $(M,F)$ is a global symmetric space, then it is an homogeneous space.
\end{corolario}
\begin{proof}
Let us consider the average metric $h$, which is necessarily globally symmetric. Then by application of the analogous Riemannian result,
$M\cong G/H$, where $G$ is the isometry group and $H$ the isotropy group of the average metric $h$.
\end{proof}
It is remarkable that our result applies also to non-reversible Finsler metrics, as our next example shows.
\begin{ejemplo}
Let us consider a Randers space \cite{Randers} of the form $F_R=\alpha+\,\beta$, with Finsler function  in local coordinates given by the expression  \cite{BCS}
\begin{align*}
F_R(x,y)=\,\sqrt{a_{ij}\,y^i\,y^j}+\,b_ i\,y^i,
\end{align*}
where $\alpha=\,\sqrt{a_{ij}\,y^i\,y^j}$ and $b_ iy^i=\beta$ is the action of the $1$-form $b$ on the tangent vector $y\in T_xM$.
 The associated fundamental tensor is
\begin{align*}
g^r_{ij}(x,y)=\,\frac{F_R}{\alpha}\,\left(a_{ij}-\frac{y_i}{\alpha}\,\frac{y_{j}}{\alpha}\right)+\,\frac{y_i}{F_R}\,\frac{y_{j}}{F_R}.
\end{align*}
On the indicatrix $I_x$ the fundamental tensor is
\begin{align}
g^r_{ij}(x,y)|_{\alpha+\beta=1}=\,\frac{1}{\alpha}\,\left(a_{ij}-\frac{y_i}{\alpha}\,\frac{y_{j}}{\alpha}\right)+\,{y_i}\,{y_{j}}.
\label{fundamental tensor on the indicatrix}
\end{align}
Formally, the right hand side of the expression \eqref{fundamental tensor on the indicatrix} does not depend upon the particular details of the $1$-form $b$. Thus we can evaluate the average and it will be equal to the Euclidean case: $b=0$, $\alpha|_{I_x}=1$. In particular, if we choose $a_{ij}=\delta_{ij}$, to simplify the argument, then we have
\begin{align}
\langle g^r_{ij}(x,y)\rangle=\,\frac{1}{vol(I_x)}\,\delta_{ij}.
\end{align}
If the space $(M,F_R)$ is Berwald, then $vol(I_x)$ is constant and the space $(M,F_R)$ is symmetric. Indeed we have that the same conclusion holds for a general space $F_R=\,\alpha+\beta$ with the function $x\mapsto vol(I)_x$ constant on $M$ and $(M,a)$ symmetric.
\begin{comentario}
This example contrasts with the main result in \cite{DengHou}, since our result shows that one can have  global symmetric spaces  which are not Berwald spaces, if the reversibility condition on the metric $F$ is not assumed. 
\end{comentario}
\begin{comentario}
Remarkably, the average metric in {\it example} \ref{Randers average metric} is independent on the $1$-form $b$, except for topology of the manifold $M$ that determines the cohomology class of the $1$-forms defined on $M$.
\end{comentario}
\label{Randers average metric}
\end{ejemplo}

\subsection{Curvature average isometric invariants}
We have briefly considered before the  average of a generic curvature endomorphism of the linear connection on $\pi^*TM$ $\nabla$. The first of these average operators
is the {\it average} $hh$-{\it curvature endomorphisms}, defined before as the endomorphism
\begin{align}
\langle R \rangle_x (X_1,X_2):& T_xM\to T_xM,\quad Y\mapsto \langle R \rangle_x (X_1 ,X_2)\,Y:=\langle R^{\nabla}(\,\iota_u({X}_1),\iota_u(X_2))\rangle Y.
\label{averageR}
\end{align}
 Let us denote the vertical lift of 
$X=X^i \frac{\partial}{\partial x ^i}\in T_x M$ by $\kappa(X)=X^i 
\frac{\partial}{\partial y ^i}\in \, \mathcal{V}_u$ with $u\in\,\hat{\pi}^{-1}(x)$. The average hv-curvature in the directions $X_1$ and $X_2$ is the endomorphism
\begin{align}
\langle P \rangle_x (X_1,X_2):& T_xM\to T_xM,\nonumber \\
 & Y\mapsto\langle P \rangle_x (X_1 ,X_2)\,Y:=\langle\pi _2 P_u (\iota _u({X_1}), \kappa _u ({X_2})) \hat{\pi}^*_u  Y\rangle,
\label{averagep}
\end{align}
with $u\in I_x\subset \hat{\pi}^{-1}(x)\subset N.$
Similarly, for the vv-curvature in the case of an arbitrary linear
connection on $\hat{\pi}^* TM$, we define the average
homomorphisms,
\begin{align}
\langle Q \rangle_x (X_1 ,X_2):& T_xM\to T_xM\nonumber \\
 & Y\mapsto \langle Q \rangle_x (X_1 ,X_2):=\langle\pi _2  Q_u (\kappa _u({X_1}), \kappa _u ({X_2}))  
\hat{\pi}^*_v  Y\rangle,
\label{averageQ}
\end{align}
with $ u\in I_x\subset \hat{\pi}^{-1}(x)\subset N.$
The average endomorphisms
\eqref{averageR} and \eqref{averageQ} are defined on the manifold $M$. In the
case of a Riemannian metric, the average curvatures $\langle P \rangle_x (X_1 ,X_2)$ and $\langle Q \rangle_x (X_1 ,X_2)$ are both zero, for any $X_1,X_2 \in \,\Gamma TM$.

 The Cartan and Chern-Rund connections are invariant under isometries of $F$, since the connections are defined in terms
of the Finsler function $F$ and the fundamental tensor $g$ (that determines the isometries). Therefore, if the measure used
in the definition of the averaging operation is invariant under  isometries of $F$, the endomorphisms
$\langle P \rangle(X_1,X_2)$ and $\langle Q \rangle(X_1,X_2)$ are also invariant under the fiber isometries.
Then for the Chern-Rund connection $Q=0$, there are  defined global affine isometric invariants of the form
\begin{align}
{Inv}(M)=\int_M d\mu \,\mathcal{F}_R(\langle R \rangle ,h),
\label{isometricinvariant1}
\end{align}
and also of the form
\begin{align}
{Inv}(M)=\int_M d\mu \,\mathcal{F}_P(\langle P \rangle ,h),
\label{isometricinvariant2}
\end{align}
where
$\mathcal{F}_R(\langle R\rangle , h)$ and $\mathcal{F}_P(\langle P \rangle, h)$ are scalar functions and
the volume form $d\mu$ is the volume form associated to the
average Riemannian  metric $\langle g \rangle$. Thus the integrals \eqref{isometricinvariant1} and \eqref{isometricinvariant2} are invariant under isometries.
\end{section}
\\
{\bf Acknowledgements}.
This work was financially supported by FAPESP  n.\ 2010/11934-6, Brazil and by PNPD-Capes n. 2265/2011, Brazil.

\footnotesize{
{}

\end{document}